\documentclass[12pt,a4paper,reqno]{amsart}
\usepackage{amsmath}
\usepackage{amssymb,latexsym}
\usepackage[dvips]{graphicx}
\usepackage{color}
\usepackage{cite}
\usepackage{amsthm}

\newtheorem{theorem}{Theorem} [section]
\newtheorem{lemma}{Lemma} [section]

\newtheorem{remark}{Remark}[section]

\let\ssection=\section\renewcommand{\section}{\setcounter{equation}{0}\ssection}
\begin{document}
\date{}

\address{M. Darwich: Lebanese University, Hadath.} \email{Mohamad.Darwich@lmpt.univ-tours.fr}

\title[Blowup]{On the $L^{2}$-critical nonlinear Schrodinger equation with an inhomogeneous  damping term.}
\author{Mohamad Darwich}
\keywords{Damped Nonlinear Schr\"odinger Equation, Blow-up, Global existence.}
\begin{abstract}
 We consider the $L^2$-critical nonlinear Schrodinger equation with an
inhomogeneous damping term.
We prove that there exists an initial data such that the corresponding solution is global in $H^1(R^d)$ and we give the minimal time of the blow up for some initial data.\\
\noindent{R\'esum\'e.}
On considere l'equation de Schrodinger nonlin\'eaire $L^2$-critique avec un terme d'amortissement non homog\`ene. On montre qu'ils existent des donn\'ees initiales tels que la solution est globale dans $H^1(\mathbb{R}^{d})$ et on donne le temps minimal d'explosion pour quelques données initiales.
\end{abstract}

\maketitle
\section {Section francaise abr\'eg\'ee}
 On montre dans cette note qu'ils existent des donn\'ees initiales tels que la solution est globale dans $H^1(\mathbb{R}^{d})$ et on donne le temps minimal d'explosion pour quelques données initiales, le resultat sera obtenue en montrant qu'un ph\'enom\`ene de $L^2$ concentration aura lieu proche de point d'explosion.\\
 Notre r\'esultat principal concernant l'equation (\ref{NLSa})est donn\'e dans le th\'eor\`eme suivant:
 \begin{theorem}
 Soit $u_0 \in H^{1}(\mathbb{R}^{d})$ et $d=1,2,3,4$:
 \begin{enumerate}
 \item Si $ a(x) >0$ et $\|u_0\|_{L^2} \leq \|Q\|_{L^2}$, alors la solution de (\ref{NLSa})  est globale dans $H^1(\mathbb{R}^{d})$.\\
  \item Si $a(x)$ est de signe quelconque et s'il existent des donn\'ees initiales $u_0$ avec $\|u_0\|_{L^2} \leq \|Q\|_{L^2}$ tels que la solution de (\ref{NLSa}) explose en temps fini $T$ alors $T > \frac{1}{\|a\|_{L^\infty}}\log(\frac{\|Q\|_{L^2}}{\|u_0\|_{L^2}})$.
  \end{enumerate}
 
 \end{theorem}
\section{Introduction}
In this paper, we study the Cauchy problem for  the $L^{2}$-critical nonlinear Schr\"{o}dinger equations:
\begin{equation}\label{NLSa}
\begin{cases}
iu_{t} + \Delta{u} +|u|^{\frac{4}{d}}u + ia(x)u =0,  (t,x) \in [0,\infty[\times \mathbb{R}^{d}, d=1,2,3,4. \\
u(0)= u_{0} \in H^1(\mathbb{R}^{d})
\end{cases}
\end{equation}
with a real inhomogeneous damping term a $\in C^1(\mathbb{R}^{d},\mathbb{R})\cap W^{1,\infty}(\mathbb{R}^{d},\mathbb{R})$ and initial data $u(0)= u_{0} \in H^1(\mathbb{R}^{d})$.Equation (\ref{NLSa}) arises in several
areas of nonlinear optics and plasma physics. The inhomogenous damping term corresponds to an
electromagnetic wave absorved by an inhomogenous medium. (cf \cite{Barontini},\cite{Brazhnyi})\\
It is known that the Cauchy problem for (\ref{NLSa}) is locally well-posed in $H^1(\mathbb{R}^{d})$(see Kato\cite{Kato} and also Cazenave\cite{Cazenave}): For any $u_{0} \in H^{1}(\mathbb{R}^{d})$, there exist $T \in (0,\infty]$ and a unique solution $u(t)$ of (\ref{NLSa}) with $u(0)=u_{0}$ such that $u \in C([0,T);H^1(\mathbb{R}^{d}))$. Moreover, T is the maximal existence time of the solution $u(t)$ in the sense that if $T < \infty$ then $\displaystyle{ \lim_{t\rightarrow T}{\|u(t)\|_{H^1(\mathbb{R}^{d})}}}=\infty$.\\
Dias and figueira \cite{Dias} studied the supercritical case($|u|^{p}u$ with $p > \frac{4}{d}$) and showed that blow-up in finite time can occur, using the virial method. In \cite{Correia}, Correia was studied the equation in dimension one, and he proved the existence of blowup phenomena in the energy space $H^{1}$\\
Let us notice that for $a=0$ (\ref{NLSa}) becomes the $L^2$-critical nonlinear Schr\"{o}dinger equation:\\
\begin{equation}\label{NLS}
\begin{cases}
 iu_{t} + \Delta u + |u|^{\frac{4}{d}}u = 0\\
 u(0)=u_{0} \in H^{1}(\mathbb{R}^{d})
 \end{cases}
 \end{equation}
 Special solutions play a fundamental role for the description of the dynamics of (\ref{NLS}). They are the solitary waves of the form $u(t, x) =\exp(it)Q(x)$, where $Q$ solves:
 \begin{equation}\label{ellip}
 \Delta Q + Q|Q|^{\frac{4}{d}} = Q.
 \end{equation}
Let $u$ be a solution of (\ref{NLSa}), we define the following quantities:\\
$L^2$norm : $\left\|u(t)\right\|_{L^2}.$\\
Energy : $E(u(t)) = \frac{1}{2}\|\nabla u\|_{L^2}^{2} - \frac{d}{4 + 2d}\|u\|_{L^{\frac{4}{d}+2}}^{\frac{4}{d}+2}.$\\
Kinetic momentum : $P(u(t))=Im(\displaystyle{\int} \nabla u \overline{u}(t,x)).$\\

 It is easy to prove that if $u$ is a solution of (\ref{NLSa}) then :
 \begin{equation}\label{mass a}
 \displaystyle{\frac{d}{dt}\|u(t)\|^2_{L^2}= -\int a(x)\vert u(t,x)\vert^2dx}, t \in [0,T),
 \end{equation}
 \begin{equation}\label{derivee de lenergie}
 \frac{d}{dt}E(u(t))=-\int a(x)\nabla u(t,x)dx + \int a(x)u(t,x)^{\frac{4}{d}+2}-\Re\int(\nabla u.\nabla a)\overline{u}dx
 \end{equation}
 and
 \begin{equation}\label{moment}
 \frac{d}{dt}P(u(t))=-2\int a(x)\Im (\nabla u\overline{u})dx, t \in [0,T).
 \end{equation}
 \begin{remark}
 Remark that if $ a(x) >0, \forall x \in \mathbb{R}^{d}$, then $\|u(t)\|_{L^2}  < \|u_0\|_{L^2}$, for all $ t \in [0,T[$.
 \end{remark}
Let us now our results:
\begin{theorem}\label{theoremessentiel}
\begin{enumerate}
Let $u_0 \in H^{1}(\mathbb{R}^{d})$ with $d=1,2,3,4$, then:
\item If $a(x) > 0$ and $\|u_0\|_{L^2} \leq \|Q\|_{L^2}$, then the corresponding solution of (\ref{NLSa}) is global in $H^1$.
\item If $a(x)$ has an arbitrary sign and if there exists an initial data $u_0 \in H^{1}$ with $\|u_0\| < \|Q\|_{L^2}$ such that the corresponding solution blows up at finite time $T$, then $T > \frac{1}{\|a\|_{\infty}}\log(\frac{\|Q\|_{L^2}}{\|u_0\|_{L^2}})$.
\end{enumerate}
\end{theorem}
\section{$L^2$-concentration }
In this section, we prove theorem \ref{theoremessentiel} by extending the proof of the $L^2$-concentration phenomen, proved by Ohta and Todorova \cite{ohta} in the radial case, to the non radial case.\\
 Hmidi and Keraani  showed in \cite{Hmidi} the  $L^2$-concentration for the equation (\ref{NLS}) without the hypothese of radiality, using the following theorem:
 \begin{theorem}\label{limsup}
 Let $(v_{n})_{n}$ be a bounded family of $H^1(\mathbb{R}^{d})$, such that:
 \begin{equation}
 \limsup_{n \rightarrow +\infty}\left\|\nabla v_{n}\right\|_{L^2(\mathbb{R}^{d})} \leq M \quad and \quad \limsup_{n \rightarrow +\infty}\left\|v_{n}\right\|_{L^{\frac{4}{d} + 2}} \geq m.
 \end{equation}
 Then, there exists $(x_{n})_{n} \subset \mathbb{R}^{d}$ such that:
 \begin{equation}
 v_{n}(\cdot + x_{n}) \rightharpoonup V \quad weakly, \nonumber
 \end{equation}
 with $\left\|V\right\|_{L^2(\mathbb{R}^{d})} \geq (\frac{d}{d+4})^{\frac{d}{4}}\frac{m^{\frac{d}{2}+1} + 1}{M^{\frac{d}{2}}}\left\|Q\right\|_{L^2(\mathbb{R}^{d})}$.
 \end{theorem}
 Now we have the following theorem:
 \begin{theorem}\label{nonradiale}
 Assume that $u_{0} \in H^{1}(\mathbb{R}^{d})$ , and suppose that the solution of (\ref{NLSa}) with $u(0)=u_{0}$ blows up in finite time $T \in (0,+\infty)$. Then, for any function $w(t)$ satisfying $w(t)\left\|\nabla u(t)\right\|_{L^2(\mathbb{R}^{d})} \rightarrow \infty$ as $t \rightarrow T$, there exists $x(t) \in \mathbb{R}^{d}$ such that, up to a subsequence,
 \begin{equation}
 \displaystyle{\limsup_{t \rightarrow T}\left\|u(t)\right\|_{L^2(\left|x - x(t)\right| < w(t))} \geq \left\|Q\right\|_{L^2(\mathbb{R}^{d})}.}\nonumber
 \end{equation}
 \end{theorem}
 To show this theorem we shall need the following lemma:
 \begin{lemma}\label{lemma ohta}
  Let $T \in (0,+\infty)$, and assume that a function $F : [0, T )\longmapsto(0,�+\infty)$ is
continuous, and $\lim_{t \rightarrow T} F(t) = +\infty$. Then, there exists a sequence $(t_{k})_{k}$ such that
$t_{k}\rightarrow T $and
\begin{equation}\label{Fsur son int}
\displaystyle{\lim_{t_k \rightarrow T}\frac{\displaystyle{\int}_{0}^{t_{k}}F(\tau)d\tau}{F(t_{k})} = 0.}
\end{equation}
\end{lemma}
For the proof see \cite{ohta}.\\

\textbf{Proof of Theorem \ref{nonradiale}}:\\
Suppose that there exist an initial data $u_0$ in $H^{1}$ with $\|u_0\| \leq \|Q\|_{L^2}$ such that the corresponding solution blows up at finite time $T$.\\
By the energy identity $(\ref{derivee de lenergie})$, we have 
\begin{equation}\label{energie}
\displaystyle{E(u(t)) = E(u_{0}) - \int_{0}^{t}H(u(\tau))d\tau, \quad t \in [0,T[.}
\end{equation}
Where $H(u(t)) =\displaystyle{ -\int a(x)\vert \nabla u(t,x)\vert^2dx + \int a(x)\vert u(t,x)\vert^{\frac{4}{d}+2}-\Re\int(\nabla u.\nabla a)\overline{u}dx}$. Let us recall the Gagliardo-Nirenberg inequality:
\begin{equation}\label{gagliardo}
\left\|u\right\|_{L^{2 + \frac{4}{d}}}^{2 + \frac{4}{d}} \leq C \|\nabla u\|_{L^2}^2\|u\|_{L^2}^{\frac{4}{d}},
\end{equation}
Note that (\ref{mass a}) gives that:
\begin{equation}\label{majorationmasse}
\|u_0\|_{L^2}e^{-\|a\|_{L^\infty}t}\leq \|u\|_{L^2} \leq \|u_0\|_{L^2}e^{\|a\|_{L^\infty}t}
\end{equation}
Now using (\ref{gagliardo}) and (\ref{majorationmasse}) we obtain that:
\begin{align}
\left|H(u(t))\right| & \leq \|a\|_{L^{\infty}}\left\|\nabla u(t)\right\|_{L^2(\mathbb{R}^{d})}^{2} + \|a\|_{L^{\infty}}\left\|u(t)\right\|_{L^{2 + \frac{4}{d}}}^{2 + \frac{4}{d}} + \|\nabla a\|_{L^\infty}\|u\|_{L^2}\|\nabla u\|_{L^2}\nonumber\\
&\leq \|a\|_{L^{\infty}}\left\|\nabla u(t)\right\|_{L^2(\mathbb{R}^{d})}^{2} + C\|a\|_{L^{\infty}}e^{\|a\|_{L^\infty}t}\|\nabla u\|_{L^2}^2\|u_0\|_{L^2}^{\frac{4}{d}} \nonumber\\&+ e^{\|a\|_{L^\infty}t}\|\nabla a\|_{L^\infty}\|u_0\|_{L^2}\|\nabla u\|_{L^2} \nonumber
\end{align}
for all $t \in [0,T[$. Then $$\left|H(u(t))\right|\leq C(\|a\|_{H^1},\|u_0\|_{L^2})e^{\|a\|_{L^\infty}t}\|\nabla u\|^2_{L^2}
$$ Moreover, we have $\displaystyle{\lim_{t \rightarrow T}\left\|\nabla u(t)\right\|_{L^2(\mathbb{R}^{d})}} = +\infty$, thus by Lemma \ref{lemma ohta}, there exists a sequence $(t_{k})_{k}$ such that $t_{k} \rightarrow T$ and 
\begin{equation}\label{ksurnabla}
\displaystyle{\lim_{k \rightarrow \infty}\frac{\displaystyle{\int}_{0}^{t_{k}}H(u(\tau))d\tau}{\left\|\nabla u(t_k)\right\|_{L^2(\mathbb{R}^{d})}^{2}} = 0.}
\end{equation}
Let
 $$\rho(t) = \frac{\left\|\nabla Q\right\|_{L^2(\mathbb{R}^{d})}}{\left\|\nabla u(t)\right\|_{L^2(\mathbb{R}^{d})}} \quad \text{and} \quad v(t,x)=\rho^{\frac{d}{2}}u(t,\rho x)$$
and $\rho_{k} = \rho(t_{k}), v_{k} = v(t_{k},.)$. The family $(v_{k})_{k}$ satisfies 

$$\left\|v_{k}\right\|_{L^2(\mathbb{R}^{d})} \leq e^{\|a\|_{L^\infty}T}\left\|u_{0}\right\|_{L^2(\mathbb{R}^{d})}\quad \text{and} \quad \left\|\nabla v_{k}\right\|_{L^2(\mathbb{R}^{d})} = \left\|\nabla Q\right\|_{L^2(\mathbb{R}^{d})}.$$
\bigskip
By (\ref{energie}) and (\ref{ksurnabla}), we have
\begin{equation}\label{Edevk}
\displaystyle{E(v_{k}) = \rho^2_{k}E(u_{0}) -\rho^2_{k}\int_{0}^{t_{k}}H(u(\tau))d\tau \rightarrow 0,}
\end{equation}
which yields 
\begin{equation}\label{vk ver Q}
\displaystyle{\left\|v_{k}\right\|_{L^{\frac{4}{d} + 2}}^{\frac{4}{d} + 2} \rightarrow \frac{d + 2}{d}\left\|\nabla Q\right\|_{L^2(\mathbb{R}^{d})}^{2}.}
\end{equation}
The family $(v_{k})_{k}$ satisfies the hypotheses of Theorem \ref{limsup} with \\
$$m^{\frac{4}{d} + 2} = \frac{d+2}{d}\left\|\nabla Q\right\|_{L^2(\mathbb{R}^{d})}^{2} \quad \text{and} \quad M = 
\left\|\nabla Q\right\|_{L^2(\mathbb{R}^{d})},$$
\bigskip
thus there exists a family $(x_{k})_{k} \subset \mathbb{R}^{d}$ and a profile $V \in H^{1}(\mathbb{R}^{d})$ with $\left\|V\right\|_{L^2(\mathbb{R}^{d})} \geq \left\|Q\right\|_{L^2(\mathbb{R}^{d})}$, such that,
\begin{equation}\label{convergencefaible}
\displaystyle{\rho^{\frac{d}{2}}_{k}u(t_{k}, \rho_{k}\cdot +  x_{k}) \rightharpoonup V \in H^{1} \quad \text{weakly}.}
\end{equation}
Using (\ref{convergencefaible}), $\forall A \geq 0$
\begin{equation}
 \displaystyle{\liminf_{n\to +\infty}\int_{B(0,A)}\rho_{n}^{d}|u(t_{n},\rho_{n}x+x_{n})|^{2}dx\geq \int_{B(0,A)}|V|^{2}dx,}\nonumber
 \end{equation}
  but $\lim_{n\to +\infty}\frac{w(t_{n})}{\rho_{n}}=+\infty$\,\,thus $\frac{w(t_{n})}{\rho_{n}}> A$, $\rho_{n}A < w(t_{n})$. This gives immediately:
  
  \begin{align}
  \displaystyle{\liminf_{n\to +\infty}\sup_{y\in\mathbb{R}^{d}}\int_{|x-y|\leq w({t_{n})}}|u(t_{n},x)|^{2}dx\geq \int_{|x|\leq A}|V|^{2}dx.}\nonumber
  \end{align}
  This it is true for all $A > 0$ thus :
  
  \begin{equation}\label{L2phenomena}
  \displaystyle{\liminf_{t\to T}\sup_{y\in\mathbb{R}^{d}}\int_{|x-y|\leq w(t)}|u(t,x)|^{2}dx\geq \int Q^{2} dx.}
  \end{equation}
  Finally we obtain:\\
  \begin{enumerate}
  
  \item If $a(x) >0$, then the norm $L^2$ is strictly decreasing, with (\ref{L2phenomena}) in hand we obtain that:$ \displaystyle{\|u_0\|_{L^2}^2 > \|Q\|_{L^2}^2}$\\
  \item If the sign of $a$ is arbitrary, (\ref{L2phenomena}) gives that $e^{2\|a\|_{L^\infty}T}\|u_0\|_{L^2}^2 \geq \|Q\|^2_{L^2}$.
  \end{enumerate}
  
  \textbf{Proof of Theorem \ref{theoremessentiel}:}\\
  Now if $a(x) >0$ $\forall x \in \mathbb{R}^{d}$, and $u_0$ be an initial data such that $\|u_0\| \leq \|Q\|$, we obtain a contradiction, that means that the solution is global in $H^{1}$, and this gives the proof of part 1 of the theorem.\\
  
 If a has an arbitrary sign, and if $u$ blows up with initial data $u_0$ with $\|u_0\|_{L^2} < \|Q\|_{L^2} $ at finite time $T$, we obtain that $T > \frac{1}{\|a\|_{L^\infty}}\log(\frac{\|Q\|_{L^2}}{\|u_0\|_{L^2}})$, this gives the proof of the second part of the theorem.

\end{document}